# A Geometric Approach to Solve Fuzzy Linear Systems

Nizami Gasilov [a], Şahin Emrah Amrahov [b], Afet Golayoğlu Fatullayev [c],
Halil Ibrahim Karakaş [d], Ömer Akın [e]

[a] Baskent University, Ankara, 06810 Turkey (phone: (+90)3122341010/1220; fax: (+90)3122341051; e-mail: gasilov@baskent.edu.tr). (corresponding author)

[b] Department of Computer Engineering, Engineering Faculty of Ankara University, Aziz Kansu Building, Tandogan Kampus, Tandogan, Ankara, 06100 Turkey (e-mail: emrah@eng.ankara.edu.tr).

[c] Baskent University, Ankara, 06810 Turkey (e-mail: afet@baskent.edu.tr).

[d] Baskent University, Ankara, 06810 Turkey (e-mail: karakas@baskent.edu.tr).

[e] TOBB ETU, Ankara, 06560 Turkey (e-mail: omerakin@etu.edu.tr).

In this paper, linear systems with a crisp real coefficient matrix and with a vector of fuzzy triangular numbers on the right-hand side are studied. A new method, which is based on the geometric representations of linear transformations, is proposed to find solutions. The method uses the fact that a vector of fuzzy triangular numbers forms a rectangular prism in $n$-dimensional space and that the image of a parallelepiped is also a parallelepiped under a linear transformation. The suggested method clarifies why in general case different approaches do not generate solutions as fuzzy numbers. It is geometrically proved that if the coefficient matrix is a generalized permutation matrix, then the solution of a fuzzy linear system (FLS) is a vector of fuzzy numbers irrespective of the vector on the right-hand side. The most important difference between this and previous papers on FLS is that the solution is sought as a fuzzy set of vectors (with real components) rather than a vector of fuzzy numbers. Each vector in the solution set solves the given FLS with a certain possibility.

The suggested method can also be applied in the case when the right-hand side is a vector of fuzzy numbers in parametric form. However, in this case, $\alpha$-cuts of the solution can not be determined by geometric similarity and additional computations are needed.

## 1. Introduction

Fuzzy linear systems arise naturally in many application problems and in their solutions, and have been studied by many researchers such as Friedman et al. [11-12], Abbasbandy et al. [1], Allahviranloo [2-3], Asady et. al. [5], Ezzaty [10], Peeva [21], Matinfar et al. [16-17], Nasseri and Khorramizadeh [20].

The concept of a fuzzy number was introduced and developed by Zadeh [23], Chang and Zadeh [6], Mizumoto and Tanaka [18], Dubois and Prade [9], Nahmias [19]. Different approaches to the concept were suggested by Puri and Ralescu [22], Goetschell and Voxman [14-15], Cong-Xin and Ming [7-8], Friedman and Kandel [13].

Friedman et al. [11] dealt with FLS with a crisp real coefficient matrix and with a vector of fuzzy numbers in parametric form on the right-hand side. They solved a $2n \times 2n$ crisp linear system obtained from a given $n \times n$ FLS. In [12], Friedman et al. studied the dual system of FLS. Further studies concentrated on solution methods. Allahviranloo [3], Matinfar et al. [16] and Matinfar et al. [17] applied the method of Adomian decomposition, the method of Householder decomposition and the Greville algorithm, respectively.

In this paper, we propose a geometric approach to FLS with a crisp real coefficient matrix and with a vector of triangular fuzzy numbers on the right-hand side. In this approach, we study the geometric object in $n$ dimensional space corresponding to the given FLS. We show that a vector of triangular fuzzy numbers represents a rectangular prism in $n$ dimensional space, and that the solutions are contained in an $n$ dimensional parallelepiped if the coefficient matrix is invertible. Then we suggest an efficient algorithm to determine this parallelepiped. We also give a method to check if a given vector is a solution and to compute its possibility if it is.

In this paper, unlike earlier research, we are not looking for solutions of FLS in the form of a vector of fuzzy numbers. In our opinion, requiring the solution to be a fuzzy number is more of a mathematical constraint rather than a natural requirement in the problem. Instead, our solutions constitute a fuzzy set of vectors of real numbers. Each vector in the solution set is a solution of the system with a certain possibility.

The suggested geometric approach also explains why solutions obtained by means of other methods are not usually vectors of fuzzy numbers. In addition, we give a geometric proof of the theorem about the necessary and sufficient condition for the existence of solution in the form of a vector of fuzzy numbers.

In fact, our approach holds for all kinds of fuzzy numbers; we choose the right-hand side as a vector of triangular fuzzy numbers to simplify computations only.

This paper is comprised of 6 sections including the introduction. Preliminaries are given in Section 2. In Section 3, we define FLS and its solution. In Section 4, we apply a geometric approach to find the solutions of FLS and present the main results. In Section 5, we solve samples of FLS. In Section 6, we summarize the results.

## 2. Preliminaries

The most popular kind of fuzzy numbers is triangular fuzzy numbers. These numbers are denoted as $\tilde{u} = (a,c,b)$ and their membership functions are defined as follows:

$$\mu(x) = \begin{cases} \dfrac{x-a}{c-a}, & a \leq x \leq c \\ \dfrac{x-b}{c-b}, & c \leq x \leq b \end{cases}$$

where $c \neq a, c \neq b$.

*Remark*: A crisp number $a$ may be regarded as the triangular fuzzy number $(a,a,a)$.

For a given possibility $\alpha \in [0,1]$, the left and the right boundaries of $\alpha$ cuts of a triangular fuzzy number $\tilde{u} = (a,c,b)$ are given by $u_L(\alpha) = a + \alpha(c-a)$ and $u_R(\alpha) = b + \alpha(c-b)$, respectively.

We denote $\underline{u} = a$ and $\bar{u} = b$.

**Definition 1.** Let $\tilde{u} = (a,c,b)$ be a triangular fuzzy number. The number $u_{cr} = c$ is called to be the crisp part of $\tilde{u}$. The uncertainty of $\tilde{u}$, denoted by $\tilde{0}$, is the triangular fuzzy number $\tilde{0} = (a-c, 0, b-c)$.

By Definition 1: $\tilde{u} = u_{cr} + \tilde{0}$.

In the definitions below, $\tilde{u} = (a,c,b)$ and $\tilde{v} = (d,f,e)$ are triangular fuzzy numbers.

**Definition 2.** $\tilde{u} = \tilde{v}$ if and only if $a = d$, $c = f$, $b = e$.

**Definition 3.** $\tilde{u} + \tilde{v} = (a+d, c+f, b+e)$

**Definition 4.** Multiplication of a triangular fuzzy number with a real number $k$ is defined as follows:
$$k\tilde{u} = \begin{cases} (ka, kc, kb), & k \geq 0 \\ (kb, kc, ka), & k < 0 \end{cases}$$

**Definition 5.** $\tilde{u} - \tilde{v} = \tilde{u} + (-1)\tilde{v}$

*Remark*: In view of Definition 3 and 5: $\tilde{u} - \tilde{v} = (a-e, c-f, b-d)$.

**Definition 6.** Let $D$ and $F$ be two fuzzy sets in $R^n$. Let $\mu_D(x)$ and $\mu_F(x)$ be the membership functions of $D$ and $F$, respectively. If $\mu_D(x) = \mu_F(x)$ whenever $x \in R^n$, then we say the fuzzy sets $D$ and $F$ are equal and we write $D = F$.

## 3. Fuzzy Linear Systems and Their Solutions

**Definition 7.** Let $a_{ij}$, $(1 \leq i, j \leq n)$ be crisp numbers and $\tilde{f}_i = (l_i, m_i, r_i)$ be triangular fuzzy numbers. The following linear system

$$\begin{cases} a_{11}x_1 + a_{12}x_2 + \ldots + a_{1n}x_n = \tilde{f}_1 \\ a_{21}x_1 + a_{22}x_2 + \ldots + a_{2n}x_n = \tilde{f}_2 \\ \vdots \\ a_{n1}x_1 + a_{n2}x_2 + \ldots + a_{nn}x_n = \tilde{f}_n \end{cases} \quad (1)$$

is called a fuzzy linear system (FLS).

One can rewrite (1) as follows using matrix notation.
$$A\tilde{X} = \tilde{B} \tag{2}$$
where $A = [a_{ij}]$ is an $n \times n$ crisp matrix and $\tilde{B} = (\tilde{f}_1, \tilde{f}_2, \ldots, \tilde{f}_n)^T$ is a vector of triangular fuzzy numbers.
We will consider systems with nonsingular matrices.

Let represent $\tilde{B}$ as sum of crisp part and uncertainty (with vertex at the origin): $\tilde{B} = \mathbf{b}_{cr} + \tilde{\mathbf{b}}$, where $\mathbf{b}_{cr} = (m_1, m_2, \ldots, m_n)^T$ and $\tilde{\mathbf{b}} = (\tilde{b}_1, \tilde{b}_2, \ldots, \tilde{b}_n)^T$ with $\tilde{b}_i = (l_i - m_i, 0, r_i - m_i) = (\underline{b}_i, 0, \overline{b}_i)$.
Note that the solution of the system (2) is in the form $\tilde{X} = \mathbf{x}_{cr} + \tilde{\mathbf{x}}$ (crisp solution + uncertainty), where $\mathbf{x}_{cr}$ is the solution of $A\mathbf{x}_{cr} = \mathbf{b}_{cr}$ and $\tilde{\mathbf{x}}$ is the solution of $A\tilde{\mathbf{x}} = \tilde{\mathbf{b}}$. The crisp solution is defined as $\mathbf{x}_{cr} = A^{-1}\mathbf{b}_{cr}$. Consequently, the matter is to determine the uncertainty $\tilde{\mathbf{x}}$.
We will find the solution of the system (2) as fuzzy set $S$ of vectors.
Let us define the following set for the vector $\tilde{\mathbf{b}}$.
$$\Pi = \{\mathbf{v} = (v_1, v_2, \ldots, v_n)^T \in R^n \mid \underline{b}_i \leq v_i \leq \overline{b}_i\}$$
This set represents an $n$ dimensional rectangular prism (with origin inside) in $R^n$.
$\mathbf{b}_{cr} + \Pi$ denotes the translation of this rectangular prism by the vector $\mathbf{b}_{cr}$. Let $\mathbf{w} \in \mathbf{b}_{cr} + \Pi$.
If $\mathbf{w} = (w_1, w_2, \ldots, w_n)^T$ then $l_i \leq w_i \leq r_i$. Each $w_i$ has an associated possibility $\alpha_i$. The number $\min_{1 \leq i \leq n} \alpha_i$ is defined to be the possibility of the vector $\mathbf{w}$ in the set $\mathbf{b}_{cr} + \Pi$.

**Definition 8.** For an invertible matrix $A$, we define the set $A^{-1}\Pi$ as $A^{-1}\Pi = \{A^{-1}\mathbf{v} \mid \mathbf{v} \in \Pi\}$.

Since the matrix $A^{-1}$ may expand, contract, reflect and/or rotate [4] the prism $\Pi$, the set $A^{-1}\Pi$ is an $n$ dimensional parallelepiped (not necessarily a rectangular prism in general).
Similarly, $\mathbf{x}_{cr} + A^{-1}\Pi$ represents a parallelepiped translated by the vector $\mathbf{x}_{cr}$.
Let $\mathbf{y} \in \mathbf{x}_{cr} + A^{-1}\Pi = A^{-1}\mathbf{b}_{cr} + A^{-1}\Pi$. Then there exists unique vector $\mathbf{w}$ such that $\mathbf{w} \in \mathbf{b}_{cr} + \Pi$ and $A^{-1}\mathbf{w} = \mathbf{y}$. The possibility of the vector $\mathbf{w}$ defines the possibility of the vector $\mathbf{y}$.

Now we define what we mean by a solution of (2).

**Definition 9.** A given vector $\mathbf{x}$ is called to be a solution of (2) with possibility $\alpha$, if corresponding unique vector $\mathbf{w} = A\mathbf{x}$ belongs to the set $\mathbf{b}_{cr} + \Pi$ with possibility $\alpha$.

**Definition 10.** We denote by $X^\alpha$ the set of all solutions with possibility $\alpha$, where $0 \leq \alpha \leq 1$. The set $S = \bigcup_{0 \leq \alpha \leq 1} X^\alpha$ is called the solution set of (2).

## 4. Main Results

**Theorem 1.** If $A$ is invertible, the solution set $S$ of (2) is $S = \mathbf{x}_{cr} + A^{-1}\Pi$.

**Proof.**

Let $\mathbf{x} \in \mathbf{x}_{cr} + A^{-1}\Pi$ with non-zero possibility $\alpha$. By definition of $A^{-1}\Pi$, there exists a $\mathbf{z} \in \Pi$ such that $\mathbf{x} = \mathbf{x}_{cr} + A^{-1}\mathbf{z} = A^{-1}(\mathbf{b}_{cr} + \mathbf{z})$. Hence, again by definition, the possibilities of the vectors $\mathbf{x}$ and $\mathbf{b}_{cr} + \mathbf{z}$ are equal. On the other hand, $\mathbf{b}_{cr} + \mathbf{z} \in \mathbf{b}_{cr} + \Pi$ and $A\mathbf{x} = \mathbf{b}_{cr} + \mathbf{z}$. Then by Definition 9, $\mathbf{x}$ is a solution with possibility $\alpha$.

Now let us take a vector $\mathbf{x}$ from the solution set $S$ with possibility $\alpha$. Then, there exists a vector $\mathbf{w}$ with possibility $\alpha$ such that $\mathbf{w} \in \mathbf{b}_{cr} + \Pi$ and $A\mathbf{x} = \mathbf{w}$. From here we obtain $\mathbf{w} = \mathbf{b}_{cr} + \mathbf{v}$ for some $\mathbf{v} \in \Pi$ (by definition of $\Pi$) and $A\mathbf{x} = \mathbf{b}_{cr} + \mathbf{v}$. That implies that $\mathbf{x} = \mathbf{x}_{cr} + A^{-1}\mathbf{v}$. Hence the vector $\mathbf{x}$ belongs to the set $\mathbf{x}_{cr} + A^{-1}\Pi$ with possibility $\alpha$.

Theorem 1 means that the solution set of (2) is $\mathbf{x}_{cr} + A^{-1}\Pi$, which is a parallelepiped translated by $\mathbf{x}_{cr}$. Possibilities of the solutions in this parallelepiped are distinct by definition. Then we need to address the following problems:

*a)* How can we compute the parallelepiped $S = \mathbf{x}_{cr} + A^{-1}\Pi$ efficiently?
*b)* How can we determine if a given vector is a solution of a given FLS? If it is, how can we find the possibility of the solution?
*c)* Is there a solution of FLS in the form of a vector of triangular fuzzy numbers? If not, why? For which matrices $A$, does such kind of solution exist?

In the rest of this Section, we work out these questions.

*a)* **An algorithm to compute the parallelepiped** $S = \mathbf{x}_{cr} + A^{-1}\Pi$

Step 1. Compute the inverse matrix $A^{-1}$.

Step 2.
  Compute the crisp vector $\mathbf{b}_{cr} = (b_1, b_2, ..., b_n)^T$ from the triangular fuzzy numbers $\tilde{f}_i = (l_i, m_i, r_i)$ with taking $b_i = m_i$. Also compute the fuzzy number vector $\tilde{\mathbf{b}} = (\tilde{b}_1, \tilde{b}_2, ..., \tilde{b}_n)^T$ by the equation $\tilde{b}_i = (l_i - m_i, 0, r_i - m_i) = (d_i, 0, h_i)$.

Step 3. Find the crisp vector $\mathbf{x}_{cr} = A^{-1}\mathbf{b}_{cr}$.

Step 4.
  Define the following vectors:

$$\mathbf{v}_1 = (d_1, 0, 0, \ldots, 0)^T \qquad \mathbf{u}_1 = (h_1, 0, 0, \ldots, 0)^T$$
$$\mathbf{v}_2 = (0, d_2, 0, \ldots, 0)^T \quad \text{and} \quad \mathbf{u}_2 = (0, h_2, 0, \ldots, 0)^T$$
$$\ldots \qquad \qquad \ldots$$
$$\mathbf{v}_n = (0, 0, \ldots, 0, d_n)^T \qquad \mathbf{u}_n = (0, 0, \ldots, 0, h_n)^T$$

These vectors determine the prism $\Pi$ completely. Let us call them support vectors of the prism. Since $h_i > 0$ and $d_i < 0$, $\mathbf{u}_i$ and $\mathbf{v}_i$ are in the same and opposite directions with the $i$-th standard basis vector $\mathbf{e}_i$, respectively.

*Remark.* Any vector in $R^n$ can be uniquely represented as a linear combination, with positive coefficients, of the support vectors.

Step 5. Constitute the fuzzy solution set:
$$S = \{\mathbf{x} = \mathbf{x}_{cr} + \sum_{i=1}^{n} \alpha_i A^{-1}\mathbf{w}_i \mid \alpha_i \in [0,1]; \mathbf{w}_i = \mathbf{v}_i \text{ or } \mathbf{w}_i = \mathbf{u}_i\} \text{ with } \mu_S(\mathbf{x}) = 1 - \max_{1 \le i \le n} \alpha_i.$$

The following theorem justifies Step 5.

**Teorem 2.** The fuzzy set $S = \{\mathbf{x} = \mathbf{x}_{cr} + \sum_{i=1}^{n} \alpha_i A^{-1}\mathbf{w}_i \mid \alpha_i \in [0,1]; \mathbf{w}_i = \mathbf{v}_i \text{ or } \mathbf{w}_i = \mathbf{u}_i\}$ with the membership function $\mu_S(\mathbf{x}) = 1 - \max_{1 \le i \le n} \alpha_i$ is the solution set of (2).

**Proof.**

Let represent $\tilde{B} = \mathbf{b}_{cr} + \tilde{\mathbf{b}}$ and $\tilde{\mathbf{b}} = (\underline{\mathbf{b}}, \mathbf{0}, \overline{\mathbf{b}})$.

We express the vectors $\underline{\mathbf{b}}$ and $\overline{\mathbf{b}}$ through standard basis vectors $\mathbf{e}_1, \mathbf{e}_2, \ldots, \mathbf{e}_n$:
$$\underline{\mathbf{b}} = \underbrace{\underline{b}_1 \mathbf{e}_1}_{\mathbf{v}_1} + \underbrace{\underline{b}_2 \mathbf{e}_2}_{\mathbf{v}_2} + \ldots + \underbrace{\underline{b}_n \mathbf{e}_n}_{\mathbf{v}_n}; \qquad \overline{\mathbf{b}} = \underbrace{\overline{b}_1 \mathbf{e}_1}_{\mathbf{u}_1} + \underbrace{\overline{b}_2 \mathbf{e}_2}_{\mathbf{u}_2} + \ldots + \underbrace{\overline{b}_n \mathbf{e}_n}_{\mathbf{u}_n}$$

Only the $i$-th coordinates of the vectors $\mathbf{v}_i$ and $\mathbf{u}_i$ are nonzero. The $i$-th coordinate of $\mathbf{v}_i$ is negative. The $i$-th coordinate of $\mathbf{u}_i$ is positive. Any crisp vector can be uniquely written as a linear combination, with positive coefficients, of the vectors $\mathbf{v}_i$ and $\mathbf{u}_i$.

The fuzzy numbers on the right-hand side of the system $A\tilde{\mathbf{x}} = \tilde{\mathbf{b}}$ constitute a rectangular prism in space:
$$\tilde{\mathbf{b}} = \Pi = \{\alpha_1 \mathbf{w}_1 + \alpha_2 \mathbf{w}_2 + \ldots + \alpha_n \mathbf{w}_n) \mid \alpha_i \in [0,1]; \mathbf{w}_i = \mathbf{v}_i \text{ or } \mathbf{w}_i = \mathbf{u}_i\}$$

This prism, upon multiplication with the inverse matrix, transforms into a parallelepiped and after translation by the vector $\mathbf{x}_{cr}$ gives the solution set:
$$S \equiv \tilde{X} = \mathbf{x}_{cr} + A^{-1}\Pi =$$
$$= \{\mathbf{x} = \mathbf{x}_{cr} + A^{-1}(\alpha_1 \mathbf{w}_1 + \alpha_2 \mathbf{w}_2 + \ldots + \alpha_n \mathbf{w}_n) \mid \alpha_i \in [0,1]; \mathbf{w}_i = \mathbf{v}_i \text{ or } \mathbf{w}_i = \mathbf{u}_i\}$$
with $\mu_X(\mathbf{x}) = 1 - \max_{1 \le i \le n} \alpha_i$.

An $\alpha$-cut of the solution, is similar to the parallelepiped above:
$$X_\alpha = \mathbf{x}_{cr} + A^{-1}\Pi_\alpha =$$
$$= \{\mathbf{x} = \mathbf{x}_{cr} + A^{-1}(\alpha_1 \mathbf{w}_1 + \alpha_2 \mathbf{w}_2 + \ldots + \alpha_n \mathbf{w}_n) \mid \alpha_i \in [0, 1-\alpha]; \mathbf{w}_i = \mathbf{v}_i \text{ or } \mathbf{w}_i = \mathbf{u}_i\}$$

Thus, the theorem has been proved.

Below we give another representation for the solution and then we extend the results for the case, where the right-hand sides of the system (2) are arbitrary fuzzy numbers.

If the triangular fuzzy numbers on the right-hand side of (2) are as $\tilde{f}_i = (l_i, m_i, r_i)$, we have $\underline{b}_i = l_i - m_i$, $\overline{b}_i = r_i - m_i$. It can be seen that an $\alpha$-cut of the solution and the solution itself can be represented also as follows:

$$X_\alpha = \{\mathbf{x} = \mathbf{x}_{cr} + c_1 A^{-1}\mathbf{e}_1 + c_2 A^{-1}\mathbf{e}_2 + \ldots + c_n A^{-1}\mathbf{e}_n \mid c_i \in [(1-\alpha)\underline{b}_i, (1-\alpha)\overline{b}_i]\};$$

$$\tilde{X} = X_0 \text{ with } \mu_X(\mathbf{x}) = 1 - \max_{1 \le i \le n} \gamma_i, \text{ where } \gamma_i = \begin{cases} c_i/\overline{b}_i, & c_i \ge 0 \\ c_i/\underline{b}_i, & c_i < 0 \end{cases}.$$

Not that $A^{-1}\mathbf{e}_i$ is determined by the $i$-th column of matrix $A^{-1}$, so no additional calculations required to constitute the parallelepipeds $X_\alpha$ and $X = X_0$.

For the case when right-hand side consists of parametric fuzzy numbers $\tilde{f}_i = ((f_i)_L(\alpha), (f_i)_R(\alpha))$ the solution can be represented as follows:

$$X_\alpha = \{\mathbf{x} = \mathbf{x}_{cr} + c_1 A^{-1}\mathbf{e}_1 + c_2 A^{-1}\mathbf{e}_2 + \ldots + c_n A^{-1}\mathbf{e}_n$$
$$\mid c_i \in [(f_i)_L(\alpha) - (b_{cr})_i, (f_i)_R(\alpha) - (b_{cr})_i]\};$$

$$\tilde{X} = X_0 \text{ with } \mu_X(\mathbf{x}) = \min_{1 \le i \le n} \alpha_i,.$$

where

$$\alpha_i = \begin{cases} (f_i)_R^{-1}(k_i), & k_i > (f_i)_R(1) \\ 1, & (f_i)_L(1) \le k_i \le (f_i)_R(1) \\ (f_i)_L^{-1}(k_i), & k_i < (f_i)_L(1) \end{cases} \text{ and } k_i = (b_{cr})_i + c_i.$$

Here $\mathbf{b}_{cr}$ is a vector with possibility 1. We note that if the initial values are in parametric form, then $\mathbf{b}_{cr}$, in general, is not unique. In this case, we can choose the components of $\mathbf{b}_{cr}$ arbitrarily to the extent that $(f_i)_L(1) \le (b_{cr})_i \le (f_i)_R(1)$. For instance, we can put $(b_{cr})_i = [(f_i)_L(1) + (f_i)_R(1)]/2$.

### b) Checking if a vector is a solution of a given FLS and computing its possibility if it is a solution.

Step 1. Compute the vector $\mathbf{z} = AX - \mathbf{b}_{cr}$.

Step 2. Represent $\mathbf{z}$ as a linear combination, with positive coefficients, of the support vectors.

$$\mathbf{z} = \sum_{i=1}^{n} \alpha_i \mathbf{w}_i, \text{ where } \mathbf{w}_i = \mathbf{v}_i \text{ or } \mathbf{w}_i = \mathbf{u}_i.$$

Step 3. If $\alpha = \max_{1 \le i \le n} \alpha_i$ is greater than 1, then $X$ is not a solution. Otherwise $X$ is a solution with possibility $1 - \alpha$.

### c) A necessary and sufficient condition for existence of a solution in the form of a vector of fuzzy numbers.

**Theorem 3.** (2) has a solution in the form of a vector of fuzzy numbers if and only if the matrix $A$ can be written as $A = DP$, where $D$ is a diagonal matrix, $P$ is a permutation matrix.

**Proof.**

In order for the solution set to correspond to a triangular fuzzy number, the set $A^{-1}\Pi$ needs to determine a rectangular prism (with edges parallel to coordinate axes). This means that either $A^{-1}$ does no rotation or it rotates by 90°, 180° or 270°. Then $A$ has the same property. These rotations can only cause the base vectors to reverse directions or replace each other. Hence, these rotations can be represented by a permutation matrix $P$. A fact in linear algebra [4] states that any linear transformation is a combination of expansion, contraction, reflection and rotation. Operations other than rotation can be represented by a diagonal matrix $D$. Hence the matrix $A$ is in the form $A = DP$.

## 5. Examples

**Example 1**: Solve the system: $\begin{bmatrix} 1 & -1 & 2 \\ 3 & -1 & 4 \\ 5 & 1 & 7 \end{bmatrix} \begin{bmatrix} x \\ y \\ z \end{bmatrix} = \begin{bmatrix} (-4, -2, -1) \\ (-1, 0, 1) \\ (12, 14, 17) \end{bmatrix}$. Determine whether

$X_1 = \begin{bmatrix} 0.5 \\ 4.5 \\ 0.9 \end{bmatrix}$ and $X_2 = \begin{bmatrix} -3.4 \\ 5.2 \\ 3.5 \end{bmatrix}$ are solutions. Compute the possibility for the solutions.

*Solution*:
We have

$$A^{-1} = \frac{1}{6}\begin{bmatrix} -11 & 9 & -2 \\ -1 & -3 & 2 \\ 8 & -6 & 2 \end{bmatrix}; \quad \mathbf{b}_{cr} = \begin{bmatrix} -2 \\ 0 \\ 14 \end{bmatrix}$$

The crisp solution is: $\mathbf{x}_{cr} = A^{-1}\mathbf{b}_{cr} = \frac{1}{6}\begin{bmatrix} -11 & 9 & -2 \\ -1 & -3 & 2 \\ 8 & -6 & 2 \end{bmatrix}\begin{bmatrix} -2 \\ 0 \\ 14 \end{bmatrix} = \begin{bmatrix} -1 \\ 5 \\ 2 \end{bmatrix}$

The uncertainty of the right-hand side is:

$\tilde{\mathbf{b}} = \tilde{B} - \mathbf{b}_{cr} = \begin{bmatrix} (-4, -2, -1) \\ (-1, 0, 1) \\ (12, 14, 17) \end{bmatrix} - \begin{bmatrix} -2 \\ 0 \\ 14 \end{bmatrix} = \begin{bmatrix} (-2, 0, 1) \\ (-1, 0, 1) \\ (-2, 0, 3) \end{bmatrix} \Rightarrow$

$\underline{\mathbf{b}} = \begin{bmatrix} -2 \\ -1 \\ -2 \end{bmatrix} = \underbrace{\begin{bmatrix} -2 \\ 0 \\ 0 \end{bmatrix}}_{\mathbf{v}_1} + \underbrace{\begin{bmatrix} 0 \\ -1 \\ 0 \end{bmatrix}}_{\mathbf{v}_2} + \underbrace{\begin{bmatrix} 0 \\ 0 \\ -2 \end{bmatrix}}_{\mathbf{v}_3}; \quad \overline{\mathbf{b}} = \begin{bmatrix} 1 \\ 1 \\ 3 \end{bmatrix} = \underbrace{\begin{bmatrix} 1 \\ 0 \\ 0 \end{bmatrix}}_{\mathbf{u}_1} + \underbrace{\begin{bmatrix} 0 \\ 1 \\ 0 \end{bmatrix}}_{\mathbf{u}_2} + \underbrace{\begin{bmatrix} 0 \\ 0 \\ 3 \end{bmatrix}}_{\mathbf{u}_3}$

$\tilde{\mathbf{b}}$ corresponds to a rectangular prism:

$\tilde{\mathbf{b}} = \Pi = \{\alpha\mathbf{w}_1 + \beta\mathbf{w}_2 + \gamma\mathbf{w}_3 \mid \alpha, \beta, \gamma \in [0, 1]; \mathbf{w}_i = \mathbf{v}_i \text{ or } \mathbf{w}_i = \mathbf{u}_i\}$.

The vectors $\mathbf{v}_1, \mathbf{v}_2, \mathbf{v}_3, \mathbf{u}_1, \mathbf{u}_2, \mathbf{u}_3$ are support vectors of the prism.

We calculate

$$\mathbf{z}_1 \equiv \mathbf{b}_1 = AX_1 - \mathbf{b}_{cr} = \begin{bmatrix} -0.2 \\ 0.6 \\ -0.7 \end{bmatrix} \text{ and } \mathbf{z}_2 \equiv \mathbf{b}_2 = AX_2 - \mathbf{b}_{cr} = \begin{bmatrix} 0.2 \\ -1.4 \\ -1.3 \end{bmatrix}$$

We express $\mathbf{b}_1$ and $\mathbf{b}_2$ as linear combinations, with positive coefficients, of support vectors:

$$\mathbf{b}_1 = \begin{bmatrix} -0.2 \\ 0.6 \\ -0.7 \end{bmatrix} = 0.1 \underbrace{\begin{bmatrix} -2 \\ 0 \\ 0 \end{bmatrix}}_{\mathbf{v}_1} + 0.6 \underbrace{\begin{bmatrix} 0 \\ 1 \\ 0 \end{bmatrix}}_{\mathbf{u}_2} + 0.35 \underbrace{\begin{bmatrix} 0 \\ 0 \\ -2 \end{bmatrix}}_{\mathbf{v}_3} \Rightarrow$$

$X_1$ is a solution with $\mu_X(X_1) = \mu_\Pi(\mathbf{b}_1) = 1 - \max\{0.1, 0.6, 0.35\} = 0.4$;

$$\mathbf{b}_2 = \begin{bmatrix} 0.2 \\ -1.4 \\ -1.3 \end{bmatrix} = 0.2 \underbrace{\begin{bmatrix} 1 \\ 0 \\ 0 \end{bmatrix}}_{\mathbf{u}_1} + 1.4 \underbrace{\begin{bmatrix} 0 \\ -1 \\ 0 \end{bmatrix}}_{\mathbf{v}_2} + 0.65 \underbrace{\begin{bmatrix} 0 \\ 0 \\ -2 \end{bmatrix}}_{\mathbf{v}_3} \Rightarrow$$

$\max\{0.2, 1.4, 0.65\} > 1 \Rightarrow X_2$ is not a solution.

The fuzzy solution set $\tilde{X}$ forms a parallelepiped in the coordinate space:

$\tilde{X} = \mathbf{x}_{cr} + A^{-1}\Pi = \{\mathbf{x} = \mathbf{x}_{cr} + A^{-1}(\alpha \mathbf{w}_1 + \beta \mathbf{w}_2 + \gamma \mathbf{w}_3) \mid \alpha, \beta, \gamma \in [0, 1]; \mathbf{w}_i = \mathbf{v}_i \text{ or } \mathbf{w}_i = \mathbf{u}_i\}$ with $\mu_X(\mathbf{x}) = 1 - \max\{\alpha, \beta, \gamma\}$

To emphasize the geometry of the problem we solve the following example using a method, equivalent to the algorithm above.

**Example 2**: Solve the system $\begin{bmatrix} 3 & 5 \\ 1 & -2 \end{bmatrix} \begin{bmatrix} x \\ y \end{bmatrix} = \begin{bmatrix} (-2, -1, 1) \\ (5, 7, 8) \end{bmatrix} = \begin{bmatrix} \tilde{a} \\ \tilde{b} \end{bmatrix}$. Find $\alpha$-cuts of the solutions for $\alpha = 0.4$ and $\alpha = 0.7$.

*Solution*:

Triangular fuzzy numbers $\tilde{a}$ and $\tilde{b}$ on the right-hand side of the system form a rectangular region in the coordinate plane. The boundary of the region shown by the rectangle *ACBD* in Fig. 1.

For $\alpha = 0.4$: $\quad a_L(0.4) = -1.6; \quad a_R(0.4) = 0.2$
$\qquad\qquad\quad b_L(0.4) = 5.8; \quad b_R(0.4) = 7.6$

For $\alpha = 0.7$: $\quad a_L(0.7) = -1.3; \quad a_R(0.7) = -0.4$
$\qquad\qquad\quad b_L(0.7) = 6.4; \quad b_R(0.7) = 7.3$

The coefficient matrix $A = \begin{bmatrix} 3 & 5 \\ 1 & -2 \end{bmatrix}$ is invertible: $A^{-1} = \frac{1}{11}\begin{bmatrix} 2 & 5 \\ 1 & -3 \end{bmatrix}$. Under the multiplication by the inverse matrix, one can find the images of the points *A*, *B*, *C* and the vertex *M*.

For this example: $A'(\frac{36}{11}, -\frac{26}{11})$, $B'(\frac{27}{11}, -\frac{14}{11})$, $C'(\frac{42}{11}, -\frac{23}{11})$ and the crisp solution is $M'(3, -2)$.

The fourth vertex of the parallelogram *A'C'B'D'* can be computed by the formula $OD' = OA' + C'B'$. Hence $D'(\frac{21}{11}, -\frac{17}{11})$.

One can determine an $\alpha$-cut of the fuzzy solution from $A', B', C'$ and $M'$ by taking $M'$, which corresponds to the crisp solution, as a center and using geometric similarity. Since the right-hand side of the system is in the form of a vector of triangular fuzzy numbers we can do this.

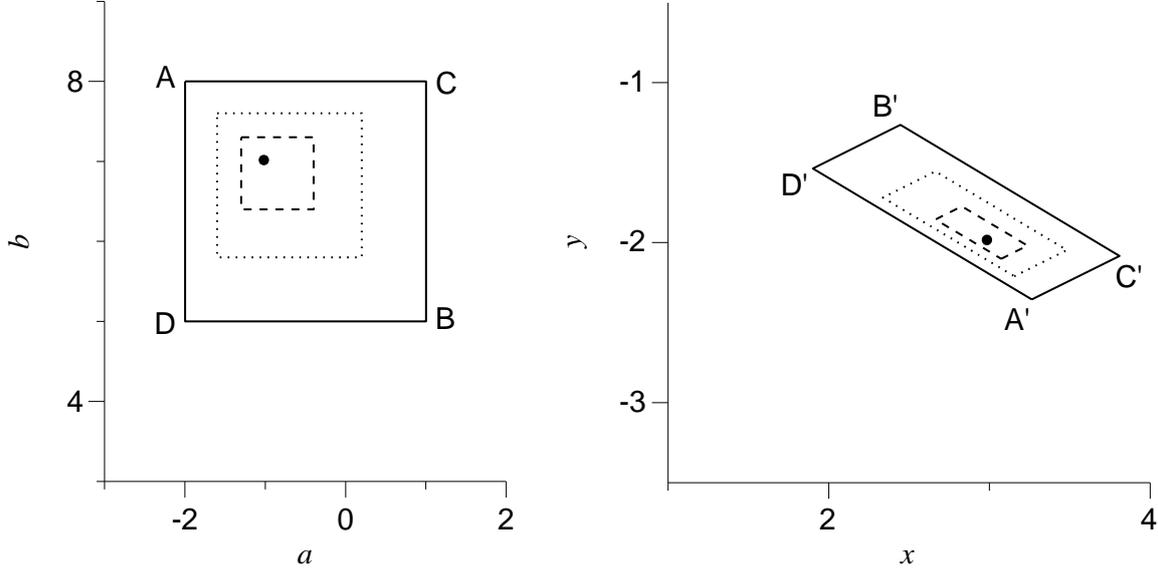

**Figure 1**. Thick black dots represent the crisp value of the right-hand side and the relevant crisp solution.
*The left part of the figure*: The rectangular boundary $ACBD$ of the fuzzy region determined by the right-hand side of the system and the boundaries of the $\alpha$-cuts ($\alpha = 0.4$ (dotted line); $\alpha = 0.7$ (dashed line)).
*The right part of the figure*: The parallelogram boundary $A'C'B'D'$ of the fuzzy region determined by the solution and the boundaries of the $\alpha$-cuts.

For $\alpha = 0.4$, we work out the vertices of the parallelogram that bounds the $\alpha$-cut.
$\hat{A} = M' + (1-\alpha)M'A' = (3, -2) + 0.6 \cdot (\frac{3}{11}, -\frac{4}{11}) = (3\frac{18}{110}, -2\frac{24}{110}) \approx (3.16364, -2.21818)$.
Similarly: $\hat{B} \approx (2.67273, -1.56364)$ and $\hat{C} \approx (3.49091, -2.05455)$. $O\hat{D} = O\hat{A} + \hat{C}\hat{B} \Rightarrow$
$\hat{D} \approx (2.34546, -1.72727)$.
For $\alpha = 0.7$:
$\breve{A} = M' + (1-\alpha)M'A' = (3\frac{9}{110}, -2\frac{12}{110}) \approx (3.08182, -2.10909)$, $\breve{B} \approx (2.83636, -1.78181)$,
$\breve{C} \approx (3.24545, -2.02727)$ and $O\breve{D} = O\breve{A} + \breve{C}\breve{B} \Rightarrow \breve{D} \approx (2.67273, -1.86363)$.

## 6. Conclusion

In this paper, we dealt with FLS with crisp coefficients and a vector of triangular fuzzy numbers on the right-hand side. We proposed a geometric approach to solve the system. Instead of looking for solutions as vectors of fuzzy numbers, we determined fuzzy set of vectors of real numbers which satisfy FLS with some possibility. We showed that this set is an *n* dimensional parallelepiped. We suggested an efficient method to compute the solution set. Then we proved the necessary and sufficient condition theorem for the existence of a solution in the form of a vector of fuzzy numbers. We illustrated the results with numerical examples.